\documentclass[11pt]{article}
\usepackage{amsmath,amssymb,amsthm,amsfonts}
\usepackage{amssymb}
\newtheorem{thm}{Theorem}[section]
\newtheorem{cor}[thm]{Corollary}
\newtheorem{lem}[thm]{Lemma}

\theoremstyle{definition}

\theoremstyle{remark}

\numberwithin{equation}{section}

\DeclareMathSymbol{\C}{\mathalpha}{AMSb}{"43}

\newcommand{\eps}{\varepsilon}

\newcommand{\lam}{\lambda}

\newcommand{\R}{{\mathbb{R}}}
\newcommand{\h}{{\mathcal{H}}}
\newcommand{\inte}{\int_{\mathbb{R}^2}}

\newcommand{\bsub}{\begin{subequations}}
\newcommand{\esub}{\end{subequations}$\!$}

\begin{document}

\title{Energy estimates and symmetry breaking  in attractive Bose-Einstein condensates
  with ring-shaped potentials 
  }

\author{\small Yujin Guo, Xiaoyu Zeng \& Huan-Song Zhou \\
 \ \\
    }
\date{}

\smallbreak \maketitle

\begin{abstract}
This paper is concerned with the properties of $L^2$-normalized minimizers of the  Gross-Pitaevskii (GP)
functional for a  two-dimensional Bose-Einstein condensate with
attractive interaction and ring-shaped potential. By establishing some delicate estimates on the least energy of the GP functional, we prove that symmetry breaking occurs for the minimizers of the GP functional as the interaction strength $a>0$ approaches a critical value $a^*$, each minimizer of the GP functional concentrates to a point on the circular  bottom of the potential well and then is non-radially symmetric as $a\nearrow a^*$. However, when $a>0$ is suitably small we prove that the minimizers of the GP functional are unique, and this unique minimizer is radially symmetric.

\end{abstract}

\vskip 0.2truein



\vskip 0.2truein

\section{Introduction}

 Since the remarkable experiments on  Bose-Einstein condensates (BEC) in dilute gases of alkali atoms in 1995 \cite{Anderson,D},  much attention has been attracted to the experimental  studies on BEC over the  last two decades, and many new phenomena of BEC have been observed  in experiments \cite{D}. These new experimental progresses also inspired  the theoretical research in BEC, especially, the theory of  Gross-Pitaevskii (GP) equations proposed  by Gross and Pitaevskii \cite{G60,G63,P}. There has been a growing interest in the mathematical theories and numerical methods of GP equations \cite{Bao2}. Several rigorous mathematical verifications of GP theory were established, see e.g. \cite{ES,LS,LSS,LSY,LSY2d}. It is known that  the classical trapping potential used in the study of BEC is the harmonic potential. With the advance of experimental techniques for BEC, some different trapping potentials have been used in the experiments \cite{CN,GM,HJ,RA,S}. Theoretically, it is also interesting to discuss mathematically  how the shapes of trapping potentials affect the behavior of BEC. Very recently, Guo and Seiringer \cite{GS} studied the BEC with attractive  interactions in $\R^2$ described by the following GP functional
 \begin{equation}
  E_a(u):=\int_{\R ^2} \big(|\nabla
  u(x)|^2+V(x)|u(x)|^2\big)dx-\frac{a}{2}\int_{\R ^2}|u(x)|^4dx\, , \quad  u\in \h \,, \label{f}
\end{equation}
 where $a>0$ describes the strength of the attractive interactions, and
\begin{equation}
   \h := \Big \{u\in  H^1(\R ^2):\ \int _{\R ^2}  V(x)|u(x)|^2 dx<\infty \Big\}\,, \label{H}
\end{equation}
with a trapping potential of the form
\begin{equation}\label{1.3}
V(x) =  h(x) \prod_{i=1}^n |x-x_i|^{p_i}, \quad p_i>0 \quad \text{and\quad $C < h(x)
< 1/C$}
\end{equation}
for some $C>0$ and all $x\in\R^2$. The authors in \cite{Bao2,GS}  proved that there exists $a^*>0$ such that the constrained minimization problem
\begin{equation}\label{eq1.3}
e(a):=\inf\Big\{E_a(u): u\in\mathcal{H} \text{ and }\int_{\mathbb{R}^2} u^2dx=1\Big\}
\end{equation}
 has at least one minimizer if and only if $a\in[0,a^*)$. Moreover,
 \begin{equation}\label{1.5}
 a^*=\inte |Q(x)|^2dx,
 \end{equation}
 and $Q(x)$ is the unique positive solution (up to translations) of the scalar field equation
  \begin{equation}
-\Delta u+ u-u^3=0\  \mbox{  in } \  \R^2,\  \  u\in
H^1(\R ^2).  \label{Kwong}
\end{equation}
The existence of $Q(x)$ is well known and $Q(x)$ is actually radially symmetric, see e.g. \cite{GNN,K,Li,mcleod}.

 In what follows, we call $e(a)$ the GP energy, which is also the least energy of a BEC system. As mentioned in \cite{GS}, the parameter $a$ in (\ref{f}) has to be interpreted as the particle number times the interaction strength, the existence of the threshold value $a^*$ described above shows that there exists a critical particle  number for collapse of the BEC \cite{D}. Theorem 1  of \cite{GS} implies that the shape of trapping potential does not affect the critical particle number. However, the behavior of the minimizers for  (\ref{eq1.3}) as $a\nearrow a^*$ does depend on the shape of potentials. In fact, for the trapping potential (\ref{1.3}), a detailed description of the behavior of the minimizers for  (\ref{eq1.3}) is given in Theorem 2 of \cite{GS}, which shows that  minimizers of (\ref{eq1.3}) must  concentrate at one of the flattest minima  $x_{i_0}(1\leq i\leq n)$ of $V(x)$ as $a\nearrow a^*$. This  also implies the presence of symmetry breaking of the minimizer. Note that the method of \cite{GS} depends heavily on the potential $V(x)$ of (\ref{1.3}) having  a finite number of minima $\{x_i\in\R^2, \ i=1,\cdots,n\}$.

 It is natural to ask what would happen if $V(x)$ has infinitely many minima.  Hence, in this paper we are mainly interested in studying the GP functional  with a trapping potential $V(x)$ with infinitely many minima and analyzing the detailed behavior of its minimizers as $a\nearrow a^*$.
 For this purpose, we focus on the following ring-shaped trapping potential:
 \begin{equation}\label{1:V}
V(x)=(|x|-A)^2,\quad \text{where}\  A>0, x\in \R^2,\end{equation}
 which is essentially an important potential used in BEC experiments, see e.g. \cite{Ha,HJ,RA}. Clearly,  all points in the  set $\{x\in\R^2:|x|=A\}$  are minima of the potential given by (\ref{1:V}).  Concerning the existence of minimizers of problem (\ref{eq1.3}), much more general potentials $V(x)$ than (\ref{1:V}) are allowed, see  \cite[Theorem 1]{GS}. But to demonstrate clearly that  symmetry  breaking does occur in the minimizers of problem (\ref{eq1.3}), the uniqueness of the minimizers of problem (\ref{eq1.3}) is used in our Corollary \ref{coro1}. So, we give first the theorem as follows.

 \begin{thm}\label{u:th}
Let $a^*$ be given by (\ref{1.5}), and let $V(x)$ be such that
\[
0\leq V(x)\in L^\infty_{\rm
loc}(\R^2)\,,\quad\lim_{|x|\to\infty} V(x) = \infty\quad\text{and}\quad\inf_{x\in
\R^2} V(x) =0\,.
\]
 Then
\begin{description}
                          \item [\bf (i)]  For all $a\in[0,a^*)$ (\ref{eq1.3}) has at least one minimizer, and there is no minimizer for (\ref{eq1.3}) if $a\geq a^*$. Moreover, $e(a) > 0$ if $a<a^*$, $\underset{a\nearrow a^*} {\lim}e(a) =
                                e(a^*) = 0$ and $e(a) = -\infty$ if $a>a^*$.
                          \item [\bf(ii)] When $a\in[0,a^*)$ is suitably  small, (\ref{eq1.3}) has a unique non-negative  minimizer in $\h$.
                        \end{description}
\end{thm}
\noindent Part (i) of the above theorem is just Theorem 1 of \cite{GS}. For part (ii), a proof based on an  implicit function theorem is given in the Appendix.

To analyze the detailed behavior of the minimizers for problem (\ref{eq1.3}), a delicate estimate on the GP functional is required. As far as we know, it is usually not easy to derive  directly the optimal energy  estimates for the GP functional (\ref{f}) under general trapping potentials. Although  the authors in \cite{GS}  developed an approach to establish  this kind of energy estimates for the potential (\ref{1.3}), it does not work well for our potential (\ref{1:V}). In fact, by following the method of \cite{GS} we are only able to get the following type of estimates
\begin{equation}\label{A:eq2.1}
C_1(a^*-a)^\frac{2}{3}<e(a)<C_2 (a^*-a)^\frac{1}{2}\quad
\mbox{as}\quad
 a\nearrow a^*,\end{equation}
see our Lemma \ref{le2.1} for the details. Therefore,  one of the aims of the paper is to provide some new ways to estimate precisely the GP energy under the potential (\ref{1:V}), which may be used effectively to handle some general type potentials. Based  on the estimates, we  may improve  the power $\frac{2}{3}$ at the left of (\ref{A:eq2.1}) to be the same as at the right, namely  $\frac{1}{2}$. That is, we have
\begin{thm}\label{thmB} Let $V(x)$ be given by (\ref{1:V}). Then, there exist two positive
constants $C_1$ and $ C_2$, independent of $a$, such that
\begin{equation}\label{energy}
C_1(a^*-a)^\frac{1}{2}<e(a)<C_2 (a^*-a)^\frac{1}{2}\quad
\mbox{as}\quad
 a\nearrow a^*.
\end{equation}
\end{thm}
With the estimates of (\ref{energy}), we may continue to analyze in detail the behavior of the minimizers of (\ref{eq1.3}) and we finally have the following theorem.

\begin{thm}\label{thmA}
Let $V(x)$ be given by (\ref{1:V}) and let $u_a$ be
a non-negative minimizer of (\ref{eq1.3}) for $a<a^*$. Given a sequence
$\{a_k\}$ with $a_k\nearrow a^*$ as $k\to\infty$, there exists a
subsequence, still denoted by $\{a_k\}$, such that each $u_{a_k}$ has a
unique maximum point $x_k$ and $x_k\to
y_0 $ as $k\to\infty$ for some $y_0\in \R^2$ satisfying  $|y_0|=A>0$.
Moreover,
\begin{equation}\label{1.10}\lim_{k\to\infty}\frac{|x_{k}|-A}{(a^*-a_k)^{\frac{1}{4}}}=
0\,,\end{equation} and
\begin{equation}
(a^*-a_k)^{\frac{1}{4}} u_{a_k}\left(x_k + x
(a^*-a_k)^{\frac{1}{4}}\right)\xrightarrow{k}  \frac { \lambda _0Q(\lambda
_0x)}{\|Q\|_2} \text{ strongly in } H^1(\R^2),\label{con:a}
\end{equation}
 where $\lam _0>0 $ satisfies
\begin{equation}\label{def:li}
  \lambda_0 = \Big(
\frac{\inte |y_0\cdot x|^2Q^2(x)dx}{A^2}\Big)^\frac{1}{4}.
\end{equation}
\end{thm}

As we mentioned before, the method used in \cite{GS} cannot be applied  directly in our case because the potential $V(x)$ in (\ref{1:V}) has infinitely many global minima and  some new difficulties must  be overcome in proving Theorem \ref{thmA}. Also,  finding the exact value of $\lambda_0$ in Theorem \ref{thmA} is more involved. Noting that the trapping potential $V(x)$ of (\ref{1:V}) is radially symmetric, it  then follows from Theorem \ref{u:th} (ii) that $e(a)$ has a unique non-negative minimizer which is also radially symmetric for small $a>0$. On the other hand, Theorem \ref{thmA} shows that any non-negative minimizer of $e(a)$ concentrates at a point on the ring $\{x\in\R^2:|x|=A\}$ as $a\nearrow a^*$, and thus it cannot be radially symmetric. This implies  that, as the strength of the interaction $a$ increases  from $0$ to $a^*$,  symmetry breaking occurs in the minimizers of $e(a)$. So we have the following corollary.

\begin{cor}\label{coro1}
Let $V(x)$ be given by  (\ref{1:V}). Then   there exist $a_*>0$ and $a_{**}>0$ satisfying  $a_{**}\leq a_{*}<a^*$ such that
\begin{description}
  \item[\bf (i)] $e(a)$ has a unique non-negative minimizer which is radially symmetric about the origin if $a\in[0,a_{**})$.
  \item [\bf (ii)] $e(a)$ has infinitely many different   non-negative minimizers, which are   not  radially symmetric if $a\in[a_*,a^*)$.
\end{description}
\end{cor}

We end this section by recalling some useful information related to the unique positive solution $Q(x)$ of (\ref{Kwong}). Take $N=2$ in (I.2) of \cite{W}, we then have  the
following Gagliardo-Nirenberg inequality:
\begin{equation}\label{GNineq}
\inte |u(x)|^4 dx\le \frac 2 {a^*} \inte |\nabla u(x) |^2dx
\inte |u(x)|^2dx ,\   \  u \in H^1(\R ^2)\,,
\end{equation}
 and  "=" holds when $u(x) = Q(x)$, where $a^*=\|Q\|_2^2$. Since $Q(x)$  is a solution of (\ref{Kwong}), it is easy to see that
\begin{equation}\label{1:id}
\inte |\nabla Q |^2dx  =\inte |Q| ^2dx=\frac{1}{2}\inte |Q| ^4dx ,
\end{equation} see also  \cite[Lemma 8.1.2]{C} for the details. Furthermore, by the results of
 \cite[Propositon 4.1] {GNN}, we know  that
 \begin{equation}
Q(x) \, , \ |\nabla Q(x)| = O(|x|^{-\frac{1}{2}}e^{-|x|}) \quad
\text{as  $|x|\to \infty$.}  \label{4:exp}
\end{equation}
Throughout the paper, we denote the norm of $L^p(\R^2)$ by $\|\cdot\|_p$ for $p\in(1,+\infty)$, and define the norms of $\h$ and $H^1(\R^2)$ by
 $$\|u\|_{\h}^2=\inte \big(|\nabla u|^2+ V(x)|u|^2\big)dx\quad\text{for } u\in\h,$$
and $$ \|u\|^2=\inte \big(|\nabla u|^2+ |u|^2\big)dx\quad\text{for } u\in H^1(\R^2), \text{ respectively.}$$
The scalar product of $\h$ is given by $$\langle u,v\rangle_{\h,\h}=\inte \nabla u\nabla v+ V(x)uvdx\quad  \text{for } u,v\in\h\,.$$

This paper is organized as follows: in Section 2 we first establish some preparatory
energy estimates and then prove  Theorem \ref{thmB}, that is,
the refined  estimates for the energy $e(a)$ are obtained.  Theorem \ref{thmA} is proved in
Section 3, where the phenomena of concentration and
symmetry breaking of the minimizers of (\ref{eq1.3}) are also discussed. Finally, by using an implicit function theorem,  Theorem \ref{u:th} (ii) is proved in the Appendix.

\section{Energy estimates: Proof of Theorem \ref{thmB}}

In this section, by using some ideas of \cite{GS}, we establish first a rough estimate like (\ref{A:eq2.1}) for the energy $e(a)$ of (\ref{eq1.3}).  Based on this estimate, some detailed properties of the minimizers of $e(a)$ can be obtained. With these we can prove Theorem  \ref{thmB}, which  gives an optimal energy estimate for $e(a)$.

\begin{lem}\label{le2.1}
Let $V(x)$ be given by  (\ref{1:V}). Then, there exist two positive constants $C_1$ and $ C_2$, independent of
$a$, such that
\begin{equation}\label{eq2.1}
C_1(a^*-a)^\frac{2}{3}\leq e(a)\leq C_2 (a^*-a)^\frac{1}{2}\quad
\mbox{as}\quad
 a\nearrow a^*.
\end{equation}
\end{lem}

\noindent{\bf Proof.} For any $\lambda>0$ and $u\in \mathcal{H}$
with $\|u\|_2^2=1$,  using  (\ref{GNineq}),
\begin{equation}\label{eq2.2}
\begin{split}
E_a(u)&\geq \int_{\mathbb{R}^2}(|x|-A)^2|
u(x)|^2dx+\frac{a^*-a}{2}\int_{\mathbb{R}^2}|u(x)|^4dx\\
&=\lambda+\int_{\mathbb{R}^2}\Big[(|x|-A)^2-\lambda \Big]|
u(x)|^2dx+\frac{a^*-a}{2}\int_{\mathbb{R}^2}|u(x)|^4dx\\
&\geq
\lambda-\frac{1}{2(a^*-a)}\int_{\mathbb{R}^2}\left[\lambda-(|x|-A)^2\right]^2_+dx,
\end{split}
\end{equation}
where $A>0$ and $[\,\cdot\,]_+ = \max\{0,\,\cdot\,\}$ denotes the
positive part. For $\lam>0$ small enough, we have
\begin{equation*}\begin{split}
 \int_{\mathbb{R}^2}\left[\lambda-(|x|-A)^2\right]^2_+dx &=
 2\pi\int^{A+\sqrt{\lambda}}_{A-\sqrt{\lambda}}[\lambda-(r-A)^2]^2rdr \\ &=
 2\pi\lambda^2\int^{\frac{\pi}{2}}_{-\frac{\pi}{2}}\cos^4\beta (A+\sqrt{\lambda}\sin \beta)\sqrt{\lambda}\cos \beta
d\beta\leq C\lambda^\frac{5}{2},
\end{split}
\end{equation*}
where we change the  variable $r=A+\sqrt{\lambda}\sin\beta$ with
$-\frac{\pi}{2}\le\beta \le \frac{\pi}{2}$ in the second identity.
The lower estimate of (\ref{eq2.1}) therefore follows from the above
estimate and (\ref{eq2.2}) by taking $\lam=[4(a^*-a)/(5C)]^{2/3}$ and $a\nearrow a^*$.

We next prove the upper estimate of  (\ref{eq2.1}) as follows.
Choose a non-negative $\varphi \in C_0^\infty(\R^2)$ such that
$\varphi(x) = 1$ for $|x|\leq 1$. For any $x_0\in \R^2$ satisfying
$|x_0|=A>0$, $\tau>0$ and $R>0$, set
\begin{equation}\label{def:trial}
u(x) =  A_{R,\tau} \frac{\tau}{\|Q\|_2} \varphi((x-x_0)/R) Q(\tau
(x-x_0)) \,,
\end{equation}
where $A_{R,\tau}>0$ is chosen so that $\inte u^2(x)dx=1$. By
scaling, $A_{R,\tau}$ depends only on the product $R\tau$, and we
have $\lim_{R\to \infty} A_{R,\tau} = 1$. In fact,
\begin{equation}\label{up1}
\frac 1{A_{R,\tau}^{2}} = \frac{1}{\|Q\|_2^2} \int_{\R^2} Q(x)^2
\varphi(x/(\tau R))^2 dx = 1 + o( (R\tau)^{-\infty}) \quad\text{as
$R\tau \to \infty$},
\end{equation}
because of the exponential decay of $Q$ in (\ref{4:exp}). Here
we use the notation $f(t) = o(t^{-\infty})$ for a function $f$
satisfying $\lim_{t\to \infty} |f(t)| t^s = 0$ for all $s>0$.

By the exponential decay of  (\ref{4:exp}), we have
\begin{equation}\label{a:eq2.5}
\begin{split}&\inte |\nabla u(x) |^2dx- \displaystyle\frac{a}{2}\inte u^4(x)dx \\
& =  \frac{\tau^2}{\|Q\|_2^2} \Big[\inte |\nabla Q(x)|^2dx- \displaystyle\frac{a}{2\|Q\|_2^2}\inte Q^4(x)dx + o((R\tau)^{-\infty} ) \Big]  \\
&=\displaystyle\frac{\tau^2}{2\|Q\|_2^2} \Big[\Big(
1-\frac{a}{\|Q\|_2^2}\Big) \inte Q^4(x)dx + o( (R\tau)^{-\infty})
\Big] \quad \mbox{as}\quad R\tau\to\infty ,
\end{split}\end{equation}
where (\ref{1:id}) is used. On the other
hand, we obtain that
\begin{equation}\label{eq2.5}
\begin{split}
\int_{\mathbb{R}^2}(|x|-A)^2| u|^2dx&\leq
\frac{1}{a^*}\int_{B_{R\tau}(0)}(|\frac{x}{\tau}+x_0|-|x_0|)^2Q^2(x)dx+ o((R\tau)^{-\infty} )\\
&\leq
C\int_{B_{R\tau}(0)}|\frac{x}{\tau}|^2Q^2(x)dx+o((R\tau)^{-\infty} )\\
&\leq \frac{C}{\tau^2}+o((R\tau)^{-\infty} )\quad \mbox{as}\quad
R\tau\to\infty .
\end{split}
\end{equation}
It then follows from (\ref{a:eq2.5}) and (\ref{eq2.5}) that
\begin{equation*}
e(a)\leq C(a^*-a)\tau^2+\frac{C}{\tau^2}+ o( (R\tau)^{-\infty})
\quad \mbox{as}\quad R\tau\to\infty .
\end{equation*}
By taking $\tau=(a^*-a)^{-\frac{1}{4}}$, the above inequality
implies the desired upper estimate of (\ref{eq2.1}). \qed

By Lemma \ref{le2.1} and using a similar procedures as the proof of  Lemma 4 of \cite{GS}, the following estimates for minimizers of $e(a)$ can be obtained. Here we omit the proof of the lemma.

\begin{lem}\label{le2.2}Let $V(x)$ be given by  (\ref{1:V}) and
suppose $u_a$ is a non-negative minimizer of (\ref{eq1.3}), then there
exists a positive constant $K$, independent of $a$, such that
\begin{equation}\label{eq2.6}
0<K (a^*-a)^{-\frac{1}{4}}\leq \int_{\mathbb{R}^2}| u_a|^4dx\leq
\frac{1}{K}(a^*-a)^{-\frac{1}{2}} \quad \text{ as }\quad a\nearrow
a^*. \qed
\end{equation}
\end{lem}



\begin{lem}\label{le2.3} For   $V(x)$ satisfying (\ref{1:V}), let  $u_a$ be a non-negative minimizer of (\ref{eq1.3}), and set
\begin{align}\label{eq2.9a}\epsilon ^{-2}_a:=\int_{\mathbb{R}^2}|\nabla u_a(x)|^2dx.
\end{align}
Then
\begin{description}
  \item [\bf(i)] $\epsilon_a \to 0$ as
$a\nearrow a^*$.
  \item [\bf(ii)] There exist a sequence $\{y_{\epsilon_a}\}\subset\R^2$
and positive constants $R_0$, $\eta$ such that the sequence
\begin{equation}\label{eq2.15}
 w_a(x):=\epsilon_a
u_a(\epsilon_a x+\epsilon_a y_{\epsilon_a})
\end{equation}
satisfies
\begin{equation}\label{eq2.16}
\liminf_{\epsilon_{a}\to0}\int_{B_{R_0}(0)}|w_a|^2dx\geq\eta>0.
\end{equation}
\item [\bf(iii)] The sequence  $\{\epsilon_a
y_{\epsilon_a}\}$ is bounded uniformly for $\epsilon_a \to 0$. Moreover, for any sequence $\{a_k\}$ with $a_k\nearrow a^*$, there exists a   convergent subsequence,
still denoted by $\{a_k\}$, such that
\begin{equation}\label{2.18}
\bar x_k:=\epsilon_{a_k}
y_{\epsilon_{a_k}}\rightarrow x_0 \quad\text{as } a_k\nearrow a^*
\end{equation}
for some $x_0\in \mathbb{R}^2$ being  a global minimum point of $V(x)$, i.e.,
$|x_0|=A>0$. Furthermore, we also have
\begin{equation}\label{2.19}
w_{a_k}\xrightarrow{k}\frac{\beta_1}{\|Q\|_2}Q(\beta_1|x-y_0|)  \text{ in } H^1(\R^2)\text{ for
some} \ y_0\in\mathbb{R}^2 \text{ and } \beta_1>0.
\end{equation}
\end{description}
\end{lem}

\noindent{\bf Proof.} {\bf(i):} Applying (\ref{GNineq}), it follows from Lemma \ref{le2.1} that
\begin{equation}\label{eq2.8}
\int_{\mathbb{R}^2}V(x)|u_a(x)|^2dx\leq e(a)\leq C_1(a^*-a)^\frac{1}{2}
\quad \text{ as }\quad a\nearrow a^*,
\end{equation}
and
\begin{equation*}
0\leq\int_{\mathbb{R}^2}|\nabla
u_a(x)|^2dx-\frac{a}{2}\int_{\mathbb{R}^2}| u_a(x)|^4dx=\epsilon
^{-2}_a-\frac{a}{2}\int_{\mathbb{R}^2}| u_a(x)|^4dx\leq
e(a)\xrightarrow{a\nearrow a^*} 0.
\end{equation*}
Then, by Lemma \ref{le2.2},
there exists a constant $m>0$, independent of $a$, such that
\begin{equation}\label{eq2.10}
0<m\epsilon ^{-2}_a\leq\int_{\mathbb{R}^2}|
u_a(x)|^4dx\leq\frac{1}{m}\epsilon ^{-2}_a\quad \text{ as }\quad
a\nearrow a^*,
\end{equation}
and   \begin{equation}\label{eq2.9}
C_2(a^*-a)^{-\frac{1}{4}}\leq\int_{\mathbb{R}^2}|\nabla
u_a(x)|^2dx\leq C_3(a^*-a)^{-\frac{1}{2}}\quad \text{ as }\quad
a\nearrow a^*.
\end{equation}
Hence, $\epsilon_a \to 0$ as
$a\nearrow a^*$, and $(i)$ is proved.

{\bf(ii):} Let
\begin{equation}
 \tilde{w}_a(x):=\epsilon_a u_a\big(\epsilon_a x\big)\label{eq2.10a}.
 \end{equation}
From (\ref{eq2.9a}) and (\ref{eq2.10}), we see  that
 \begin{equation}\label{eq2.12}
\int_{\mathbb{R}^2}|\nabla \tilde{w}_a|^2dx=1,\quad m\leq
\int_{\mathbb{R}^2}| \tilde{w}_a|^4dx\leq \frac{1}{m}.
\end{equation}

We claim that \begin{equation}\label{eq2.14}
\liminf_{\epsilon_a\to
0}\int_{B_{R_0}(y_{\epsilon_a})}|\tilde{w}_a|^2dx\geq\eta>0.
\end{equation}
On the contrary, suppose (\ref{eq2.14}) is false. Then for any
$R>0$, there exists a sequence $\{\tilde{w}_{a_k}\}$ with $a_k\nearrow
a^*$ such that
\begin{equation*}
\lim_{k\rightarrow\infty}\sup_{y\in
\mathbb{R}^2}\int_{B_{R}(y)}|\tilde{w}_{a_k}|^2dx=0.
\end{equation*}
By Lemma I.1 in \cite{l2} or Theorem 8.10 in \cite{Lieb}, we then
deduce that $\tilde{w}_{a_k}\xrightarrow{k}0$ in $L^p(\mathbb{R}^2)$
for any $2<p<\infty$, which  contradicts
(\ref{eq2.12}). Thus  (\ref{eq2.14}) holds,  and then (\ref{eq2.16}) follows from (\ref{eq2.10a}) and (\ref{eq2.14}).

{\bf(iii):} By (\ref{eq2.8}), we see  that
\begin{equation}\label{eq2.17}
\int_{\mathbb{R}^2}V(x)|u_a(x)|^2dx=\int_{\mathbb{R}^2}V(\epsilon_a
x+\epsilon_a y_{\epsilon_a})|w_a(x)|^2dx\rightarrow 0 \quad \text{as}\quad
a\nearrow a^*.
\end{equation}
Suppose $\{\epsilon_a y_{\epsilon_a}\}$ is unbounded as
$\epsilon_a \to 0$. Then there exists a subsequence $\{a_n\}$, where
$a_n\nearrow a^*$ as $n\to\infty $, such that
$$\epsilon_n:=\epsilon_{a_n} \to 0 \quad \text{and}\quad \epsilon_n
y_{\epsilon_n}\longrightarrow\infty \quad \text{as}\quad n\to\infty
.$$ Since $V(x)\rightarrow\infty$ as $|x|\rightarrow\infty$, we derive from (\ref{eq2.16}) and
Fatou's Lemma that for any $C>0$,
\begin{equation*}
\lim_{n\rightarrow\infty}\int_{\mathbb{R}^2}V(\epsilon_n
x+\epsilon_n y_{\epsilon_n})|w_{a_n}(x)|^2dx\geq
\int_{\mathbb{R}^2}\lim_{n\rightarrow\infty}V(\epsilon_n
x+\epsilon_n y_{\epsilon_n})|w_{a_n}(x)|^2dx\geq C \eta>0,
\end{equation*}
which contradicts (\ref{eq2.17}). Thus, $\{\epsilon_a y_{\epsilon_a}\}$ is
bounded uniformly for $\epsilon_a \to 0$. Therefore,  for any
sequence $\{a_k\}$  there exists a   convergent subsequence, still
denoted by $\{a_k\}$, such that
$$\bar x_k :=\epsilon_ky_{\epsilon_k}\xrightarrow{k} x_0\text{ for some point $x_0\in \R^2$, where $\epsilon_k:=\epsilon_{a_k}$.}$$
We claim that $|x_0|= A$. Otherwise,  if $|x_0|\neq A$, then
 $V(x_0)>0$. By applying Fatou's Lemma again, we deduce from (\ref{eq2.16}) that
\begin{equation*}
\lim_{k\rightarrow\infty}\int_{\mathbb{R}^2}V(\epsilon_k
x+\epsilon_k y_{\epsilon_k})|w_{a_k}(x)|^2dx\geq
\int_{\mathbb{R}^2}\lim_{k\rightarrow\infty}V(\epsilon_k
x+\epsilon_k y_{\epsilon_k})|w_a(x)|^2dx\geq \frac{\eta}{2}V(x_0)
>0,
\end{equation*}
which contradicts (\ref{eq2.17}). So,  $|x_0|=A>0$, and (\ref{2.18}) is proved.

We now  turn  to proving (\ref{2.19}).
Since $u_a$ is a non-negative minimizer of (\ref{eq1.3}), it satisfies
the Euler-Lagrange equation
\begin{equation}\label{eq2.18}
-\Delta u_a(x)+V(x)u_a(x)=\mu_a u_a(x)+a u_a^3(x)\quad  \text{ in }\
\mathbb{R}^2,
\end{equation}
where $\mu_a\in \mathbb{R}$ is a Lagrange multiplier, and
\begin{equation*}
\mu_a=e(a)-\frac{a}{2}\int_{\mathbb{R}^2}|u_a|^4dx.
\end{equation*}
It then follows from Lemma \ref{le2.1} and (\ref{eq2.10}) that there
exist two positive constants $C_1$ and $C_2$, independent of $a$,
such that
$$-C_2<\epsilon_a^2\mu_a <-C_1<0\quad \text{as}\quad a\nearrow a^*.$$
In view of (\ref{eq2.18}), $w_a(x)$ defined in (\ref{eq2.15})
satisfies the elliptic equation
\begin{equation}\label{eq2.19}
-\Delta w_a(x)+\epsilon_a^2V(\epsilon_a x+\epsilon_a
y_{\epsilon_a})w_a(x)=\epsilon_a^2\mu_a w_a(x)+a w_a^3(x)\quad \text { in
} \ \mathbb{R}^2.
\end{equation}
Therefore, for the convergent subsequence $\{a_k\}$  obtained in (\ref{2.18}),  we may assume
that $\epsilon_k^2\mu_{a_k}\xrightarrow{k}-\beta _1^2 < 0$ for some $\beta _1>0$, and  $w_{a_k}\overset k\rightharpoonup
w_0\geq 0$  weakly in $H^1(\mathbb{R}^2)$ for some $w_0\in H^1(\R^2)$. Since $\{\epsilon_{a}y_{\epsilon_a}\}$ is bounded uniformly in $\epsilon_{a}$, by
passing to the weak limit of (\ref{eq2.19}),
we see that   $w_0\geq0$ satisfies
\begin{equation}\label{eq2.20}
-\Delta w_0(x)=-\beta_1^2 w_0(x)+a^* w_0^{3}(x)\quad \text { \ in }
\ \mathbb{R}^2.
\end{equation}
Furthermore, it follows from (\ref{eq2.16}) that $w_0\not\equiv 0$,
and therefore we have  $w_0>0$ by the strong maximum principle. By a
simple rescaling,  the uniqueness (up to translations) of positive
solutions for the nonlinear scalar field equation (\ref{Kwong})
implies that
\begin{equation}\label{eq2.21}
w_0(x)=\frac{\beta_1}{\|Q\|_2}Q(\beta_1|x-y_0|) \quad \text{for
some} \quad y_0\in\mathbb{R}^2,
\end{equation}
where $||w_0||_2^2=1$. By the norm preservation we further conclude
that $w_{a_k}$ converges to $w_0$ strongly in $L^2(\mathbb{R}^2)$ and in
fact, strongly in $L^p(\mathbb{R}^2)$ for any $2\leq p<\infty$
because of $H^1(\mathbb{R}^2)$ boundedness. Also, since $w_{a_k}$ and
$w_0$ satisfy (\ref{eq2.19}) and (\ref{eq2.20}), respectively, a
simple analysis shows that $w_{a_k}$ converges to $w_0$
strongly in $H^1(\mathbb{R}^2)$, and thus (\ref{2.19}) holds.
\qed\\

\begin{lem}\label{le2.4}  Under the assumptions of Lemma \ref{le2.3}, let
$\{a_k\}$ be the    convergent subsequence given by Lemma \ref{le2.3} (ii). Then, for any $R>0$, there exists $C_0(R)>0$,
independent of $a_k$, such that
\begin{equation}\label{eq2.25}
\lim_{\epsilon_{a_k}\rightarrow0}\frac{1}{\epsilon_{a_k}^2}\int_{B_R(0)}V(\epsilon_{a_k}
x+\epsilon_{a_k} y_{\epsilon_{a_k}})|w_{a_k}(x)|^2dx \geq C_0(R).
\end{equation}

\end{lem}

\noindent{\bf Proof.}
Since $V(x)=(|x|-A)^2$ with $A>0$, we have
\begin{equation}\label{eq2.26}
V(\epsilon_a x+\epsilon_a
y_{\epsilon_a})=\epsilon_a^2\Big(|x+y_{\epsilon_a}|-\frac{A}{\epsilon_a}\Big)^2,
\end{equation}
where the term $|x+y_{\epsilon_a}|$  can be rewritten as
\begin{equation}\label{eq2.27}
|x+y_{\epsilon_a}|= \sqrt{|x|^2+|y_{\epsilon_a}|^2}\sqrt{1+\frac{2x\cdot
y_{\epsilon_a}}{|x|^2+|y_{\epsilon_a}|^2}}.
\end{equation}
For the convergent sequence $\{a_k\}$ given by Lemma \ref{le2.3} (ii), since $\epsilon_k y_{\epsilon_k}\xrightarrow{k} x_0$ with $|x_0|=A>0$, we have $|y_{\epsilon_k} |\xrightarrow{k}\infty$, and hence $\frac{2x\cdot
y_{\epsilon_k}}{|x|^2+|y_{\epsilon_k}|^2}\xrightarrow{k}0$ uniformly for $x\in
B_R(0)$. Using the Taylor expansion we
obtain that
$$\sqrt{1+\frac{2x\cdot y_{\epsilon_k}}{|x|^2+|y_{\epsilon_k}|^2}}=1+\frac{x\cdot y_{\epsilon_k}}{|x|^2+|y_{\epsilon_k}|^2}+O\big(\frac{1}{|y_{\epsilon_k}|^2}\big)\quad \text{for all}\quad x\in B_R(0),$$
which, together with (\ref{eq2.26}) and (\ref{eq2.27}),  then
implies that
\begin{equation}\label{eq2.28}
\frac{1}{\epsilon_k^2}V(\epsilon_k x+\epsilon_k
y_{\epsilon_k})=\Big|\sqrt{|x|^2+|y_{\epsilon_k}|^2}+\frac{x\cdot
y_{\epsilon_k}} {\sqrt{|x|^2+|y_{\epsilon_k}|^2}}-\frac{A}{{\epsilon_k}}+
O\big(\frac{1}{|y_{\epsilon_k}|}\big)\Big|^2.
\end{equation}
For any $x\in \mathbb{R}^2$, let {\it arg $x$} be the angle between
$x$ and the positive {\it x}-axis, and $\langle x,y\rangle$   the angle
between the vectors $x$ and $y$. Without loss of generality, we may assume that
$x_0=(A,0)$, and it then follows that $\arg y_{\epsilon_k}\xrightarrow{k}0$. Thus, we  can choose   $0<\delta<\frac{\pi}{16}$ small enough such that
\begin{equation}\label{eq2.29}
-\delta<\arg y_{\epsilon_k}<\delta  \quad \text{as}\quad
\epsilon_k\rightarrow0.
\end{equation}
Denote
\begin{equation}
\begin{split}
\Omega_{\epsilon_k}^1=&\Big\{x\in B_R(0):\,
\sqrt{|x|^2+|y_{\epsilon_k}|^2}\leq\frac{A}{\epsilon_k}\Big\}\\=&\Big\{x\in
B_R(0):\,|x|^2\leq\Big(\frac{A}{\epsilon_k}\Big)^2-|y_{\epsilon_k}|^2\Big\},
\end{split}
\end{equation}
and
\begin{equation}
\begin{split}
\Omega_{\epsilon_k}^2=&\Big\{x\in
B_R(0):\,\sqrt{|x|^2+|y_{\epsilon_k}|^2}>\frac{A}{\epsilon_k}\Big\}\\
=&\Big\{x\in
B_R(0):\,\Big(\frac{A}{\epsilon_k}\Big)^2-|y_{\epsilon_k}|^2<|x|^2<R^2\Big\},
\end{split}
\end{equation}
so that $B_R(0)=\Omega_{\epsilon_k}^1\cup \Omega_{\epsilon_k}^2$ and
$\Omega_{\epsilon_k}^1\cap\Omega_{\epsilon_k}^2=\emptyset$. Since
$$|\Omega_{\epsilon_k}^1|+|\Omega_{\epsilon_k}^2|=|B_R(0)|=\pi R^2,$$
there exists a subsequence, still denoted by $\{\epsilon_k\}$, of $\{\epsilon_k\}$
 such that$$\text{either}\quad
|\Omega_{\epsilon_k}^1|\geq\frac{\pi R^2}{2}\quad \text{or}\quad
|\Omega_{\epsilon_k}^2|\geq\frac{\pi R^2}{2}.$$
We finish the proof by considering the following two cases:

 {\em Case 1:} $|\Omega_{\epsilon_k}^1|\geq\frac{\pi R^2}{2}.$
In this case, we have $B_{\frac{R}{\sqrt{2}}}(0)\subset
\Omega_{\epsilon_k}^1$, and set
$$\Omega_1:=\Big(B_{\frac{R}{\sqrt{2}}}(0)\setminus B_{\frac{R}{2}}(0)\Big)\cap\left\{x:
\frac{\pi}{2}+2\delta<\arg x<\frac{3\pi}{2}-2\delta\right\}\subset
\Omega_{\epsilon_k}^1.$$ Then
\begin{equation}
|\Omega_1|=\frac{(\pi-4\delta)}{8}R^2.
\end{equation}
By (\ref{eq2.29}), one can easily check that for any $x\in
\Omega_1$,
\begin{equation}
x\cdot y_{\epsilon_k}=|x||y_{\epsilon_k}|\cos\langle
x,y_{\epsilon_k}\rangle<0\quad \text{ and }\quad |\cos\langle
x,y_{\epsilon_k}\rangle|>-\cos(\frac{\pi}{2}+\delta)>0.
\end{equation}
We thus derive from (\ref{eq2.28}) that
\begin{equation*}
\begin{split}
&\sqrt{|x|^2+|y_{\epsilon_k}|^2}+\frac{x\cdot
y_{\epsilon_k}}{\sqrt{|x|^2+|y_{\epsilon_k}|^2}}-\frac{A}{\epsilon_k}+O\big(\frac{1}{
|y_{\epsilon_k}|}\big)
\leq
\frac{x\cdot y_{\epsilon_k}}{\sqrt{|x|^2+|y_{\epsilon_k}|^2}}+O\big(\frac{1}{|y_{\epsilon_k}|}\big)\\
&\leq \frac{x\cdot
y_{\epsilon_k}}{2\sqrt{|x|^2+|y_{\epsilon_k}|^2}}\leq\frac{|x||y_{\epsilon_k}|\cos(\frac{\pi}{2}+\delta)}{2\sqrt{|x|^2+|y_{\epsilon_k}|^2}}<0\quad
\text{for }\quad x\in \Omega_1.
\end{split}
\end{equation*}
Noting that $\underset{\epsilon_k\to0}{\lim}|y_{\epsilon_k}|=\infty$,
we thus have
\begin{equation}
\frac{1}{\epsilon_k^2}V(\epsilon_k x+\epsilon y_{\epsilon_k})\geq
\frac{\cos^2(\frac{\pi}{2}+\delta)|x|^2}{8}\quad\text{for}\quad x\in
\Omega_1.
\end{equation}
Taking $\delta =\frac{\pi}{20}$, the above estimate implies that
\begin{equation}\label{eq2.35}
\begin{split}
\lim_{\epsilon_k\rightarrow0}&\frac{1}{\epsilon_k^2}\int_{B_R}V(\epsilon_k
x+\epsilon_k y_{\epsilon_k})|w_a(x)|^2dx\geq
\lim_{\epsilon_k\rightarrow0}\frac{1}{\epsilon_k^2}\int_{\Omega_1}V(\epsilon_k
x+\epsilon_k y_{\epsilon_k})|w_a(x)|^2dx\\
&\geq
\frac{\cos^2\frac{11\pi}{20}}{8}\int_{\Omega_1}|x|^2|w_0(x)|^2dx:=C(R)>0,
\end{split}
\end{equation}
and then (\ref{eq2.25}) is proved.

 {\em Case 2:} $|\Omega_{\epsilon_k}^2|\geq\frac{\pi R^2}{2}.$ In this case, we  deduce that the annular region $D_R:=B_R\setminus
B_{\frac{R}{\sqrt{2}}}\subset \Omega_{\epsilon_k}^2$. Set
$\Omega_2:=D_R\cap\left\{x; -\frac{\pi}{2}+2\delta<\arg
x<\frac{\pi}{2}-2\delta\right\}\subset \Omega_{\epsilon_k}^2$. Then,
\begin{equation}
|\Omega_2|=\frac{(\pi-4\delta)}{4}R^2.
\end{equation}
One can  check that for any $x\in \Omega_2$,
\begin{equation}\label{eq2.37}
x\cdot y_{\epsilon_k}=|x||y_{\epsilon_k}|\cos\langle
x,y_{\epsilon_k}\rangle>0\quad \text{ and }\quad \cos\langle
x,y_{\epsilon_k}\rangle>\cos(\frac{\pi}{2}-\delta)>0.
\end{equation}
It then follows from (\ref{eq2.28}) and  (\ref{eq2.37}) that
\begin{equation*}
\begin{split}
&\sqrt{|x|^2+|y_{\epsilon_k}|^2}+\frac{x\cdot
y_{\epsilon_k}}{\sqrt{|x|^2+|y_{\epsilon_k}|^2}}-\frac{A}{\epsilon_k}
+O\big(\frac{1}{|y_{\epsilon_k}|}\big)
\geq
\frac{x\cdot y_{\epsilon_k}}{\sqrt{|x|^2+|y_{\epsilon_k}|^2}}+O\big(\frac{1}{|y_{\epsilon_k}|}\big)\\
&\geq \frac{x\cdot
y_{\epsilon_k}}{2\sqrt{|x|^2+|y_{\epsilon_k}|^2}}\geq\frac{|x||y_{\epsilon_k}|\cos(\frac{\pi}{2}-\delta)}{2\sqrt{|x|^2+|y_{\epsilon_k}|^2}}>0\quad
\text{for}\quad x\in \Omega_2.
\end{split}
\end{equation*}
Hence
\begin{equation*}
\frac{1}{\epsilon_k^2}V(\epsilon_k x+\epsilon_k y_{\epsilon_k})\geq
\frac{\cos^2(\frac{\pi}{2}-\delta)|x|^2}{8}\quad \text{for }\quad
x\in \Omega_2.
\end{equation*}
Thus, by taking $\delta =\frac{\pi}{20}$, the above estimate gives
that
\begin{equation}\label{eq2.39}
\begin{split}
\lim_{\epsilon_k\rightarrow0}&\frac{1}{\epsilon_k^2}\int_{B_R}V(\epsilon_k
x+\epsilon_k y_{\epsilon_k})|w_a(x)|^2dx\geq
\lim_{\epsilon_k\rightarrow0}\frac{1}{\epsilon_k^2}\int_{\Omega_2}V(\epsilon_k
x+\epsilon_k y_{\epsilon_k})|w_a(x)|^2dx\\
&\geq
\frac{\cos^2\frac{9\pi}{20}}{8}\int_{\Omega_2}|x|^2|w_0(x)|^2dx:=C_0(R)>0.
\end{split}
\end{equation}
Therefore, (\ref{eq2.25}) also follows from (\ref{eq2.35}) and (\ref{eq2.39}) in this case.
\qed\\

We end this section by  proving  Theorem \ref{thmB}, that is, establishing  the refined estimates for $e(a)$.

\noindent \textbf{Proof of Theorem \ref{thmB}:}
By
Lemma \ref{le2.1}, it suffices to prove that
 there exists a positive $C>0$, independent of $a$, such that
\begin{equation}\label{eq2.24}
e(a)>C(a^*-a)^\frac{1}{2} \quad\text{as}\quad a\nearrow a^*.
\end{equation}
Actually, from the proof of Lemma \ref{le2.3} (iii), we see that for any sequence $\{a_k\}$ with $a_k\nearrow a^*$, there exists a   convergent subsequence, still denoted by $\{a_k\}$, such that
$w_{a_k}\rightarrow w_0>0$ strongly in
$L^4(\mathbb{R}^2)$, where $w_0$ satisfies (\ref{eq2.21}).  This implies that   there
exists a constant $M_1>0$, independent of $a_k$, such that
 $$\int_{\mathbb{R}^2}|w_{a_k}(x)|^4dx>M_1 \quad\text{as}\quad a_k\nearrow a^*.$$
Moreover, applying (\ref{eq2.25}) with $R=1$ yields that there
exists a constant $M_2>0$, independent of $a_k$, such that
\[
\int_{B_1(0)}V(\epsilon_k x+\epsilon_k y_{\epsilon_k})|w_{a_k}(x)|^2dx \geq
M_2\epsilon_k^2 \quad\text{as}\quad a_k\nearrow a^* .\] Thus,
 \begin{equation}\label{eq2.255}
 \begin{split}
e(a_k) &=E_{a_k}(u_{a_k}) =\frac{1}{\epsilon
_k^2}\left[\int_{\mathbb{R}^2}|\nabla
w_{a_k}(x)|^2dx-\frac{a^*}{2}\int_{\mathbb{R}^2}|w_{a_k}(x)|^4dx\right]\\
&+\frac{a^*-{a_k}}{2\epsilon
_k^2}\int_{\mathbb{R}^2}|w_{a_k}(x)|^4dx+
\int_{\mathbb{R}^2}V(\epsilon _k
x+\epsilon _ky_{\epsilon_k})|w_{a_k}(x)|^2dx\\
&\geq  \frac{a^*-{a_k}}{2\epsilon _k^2}M_1+M_2\epsilon _k^2 \geq
\sqrt{2M_1M_2}(a^*-{a_k})^\frac{1}{2} \quad\text{as}\quad
{a_k}\nearrow a^* ,
\end{split}
\end{equation}
and  (\ref{eq2.24}) therefore holds for the   subsequence  $\{a_k\}$.

We finally claim  that (\ref{eq2.24}) actually holds for any $a\nearrow
a^*$. Otherwise,  there exists a
sequence $\{a_m\}$ such that $a_m\nearrow a^*$ as $m\to \infty$ and
\begin{equation}\label{eq2.255A}\lim_{m\to \infty}\frac{e(a_m)}{(a^*-a_m)^{\frac{1}{2}}}=0 \,.\end{equation}
By the proof of Lemma \ref{le2.3} (iii), we know  that
there exists  a subsequence, still denoted by  $\{a_m\}$, such that $\bar x_m:=\epsilon_{m}y_{\epsilon_{m}}\xrightarrow{m}
\bar x_0$ with $|\bar x_0|=A>0$ and
$w_m(x):=w_{a_m}(x)$ defined in (\ref{eq2.15}) satisfies
\begin{equation}\label{2:wm}
w_m(x)\xrightarrow{m} \tilde{w_0}(x)>0 \quad \text{strongly in}\quad H^1(\R^2),
\end{equation}
where $\tilde{w_0}>0$ is a  solution of (\ref{eq2.20}) for some
$\beta_1>0$. Moreover, by  the same argument of  Lemma \ref{le2.4}, we also have
\begin{equation}\label{2:lim}
\lim_{m\to
\infty}\frac{1}{\epsilon^2_m}\int_{B_R(0)}V(\epsilon_m
x+\epsilon_m y_{\epsilon_m})|w_m(x)|^2dx \geq C_0(R),
\end{equation}
where $C_0(R)>0$ is independent of $m$. Applying (\ref{2:wm}) and
(\ref{2:lim}), an argument similar to (\ref{eq2.255}) yields that
there exists some constant $C>0$, independent of $m$, so that
$$e(a_m)\geq C (a^*-a_m)^\frac{1}{2} \quad\text{as}\quad
{a_m}\nearrow a^*,$$ which however contradicts (\ref{eq2.255A}).
This completes the proof of  Theorem \ref{thmB}. \qed

\section{Mass concentration and symmetry breaking}

In this section we complete the proof of Theorem \ref{thmA} under
the ring-shaped potential $V(x)=(|x|-A)^2$ with $A>0$, which addresses the mass concentration and symmetry breaking of minimizers as $a\nearrow a^*$. Let $u_a$
 be a non-negative minimizer of (\ref{eq1.3}). Applying the  energy
estimates of Theorem  \ref{thmB}, a proof  similar to that of Lemma 4 in
\cite{GS} yields that there exists a positive constant $M$,
independent of $a$, such that
\begin{equation}\label{3:1}
0<M (a^*-a)^{-\frac{1}{2}}\leq \int_{\mathbb{R}^2}| u_a|^4dx\leq
\frac{1}{M}(a^*-a)^{-\frac{1}{2}} \quad \text{ as }\quad a\nearrow
a^*.
\end{equation}
Stimulated by above estimates, we define
\begin{equation}\label{3:eps}
 \eps_a :=(a^*-a)^{\frac{1}{4 }}>0 \,.
\end{equation}
From (\ref{GNineq}) we conclude that
\[
e(a)\ge \Big(1-\frac{a}{a^*}\Big)\inte |\nabla u_a(x)|^2dx+\inte
(|x|-A)^2u^2_a(x) dx\,,
\]
and it hence follows from Theorem \ref{thmB} that
\begin{equation}
 \inte |\nabla u_a(x)|^2dx \le C\eps_a ^{-2}\quad \mbox{and}\quad \inte (|x|-A)^2u_a^2(x)dx \le C\eps_a ^2.
\label{3:2}
\end{equation}
Similar to Lemma \ref{le2.3} (ii), for $\eps_a$ given by (\ref{3:eps}), we know  that there
exist a sequence $\{y_{\varepsilon_a}\}\subset\R^2$ and positive constants $R_0$
and $\eta$ such that
\begin{equation}\label{3:eq2.16}
\liminf_{\varepsilon_a\searrow
0}\int_{B_{R_0}(0)}|w_a|^2dx\geq\eta>0,
\end{equation}
where we define the $L^2(\mathbb{R}^2)$-normalized
function\begin{equation}\label{3:eq2.15}
 w_a(x)=\varepsilon_a
u_a(\varepsilon_a x+\varepsilon_a y_{\varepsilon_a}).
\end{equation}
Note from (\ref{3:1}) and (\ref{3:2}) that
 \begin{equation}\label{3:4}
\int_{\mathbb{R}^2}|\nabla {w}_a|^2dx<C,\quad M\leq
\int_{\mathbb{R}^2}| {w}_a|^4dx\leq \frac{1}{M},
\end{equation}
where the positive constants $C$ and $M$ are independent of $a$.

\begin{lem}\label{le3.1}
For any given sequence $\{a_k\}$ with $a_k\nearrow a^*$, let $\varepsilon
_k:=\varepsilon_{a_k}=(a^*-a_k)^{\frac{1}{4 }}>0$, $u_k(x):=u_{a_k}(x)$  be a non-negative
minimizer   of (\ref{eq1.3}), and $w_k:=w_{a_k}\ge 0$ be defined by (\ref{3:eq2.15}). Then, there is a subsequence, still denoted by $\{a_k\}$, such that
\begin{equation}\label{3.7}z_k:=\varepsilon _ky_{\varepsilon _k}\xrightarrow{k}y_0 \text{ for some $y_0\in\R^2$ and  $|y_0|=A$ }.\end{equation}
Moreover, for any $\delta>0$ small enough, we have
\begin{equation}\label{3.8}
u_k(x)=\frac{1}{\varepsilon_k}w_k(\frac{x-z_k}{\varepsilon_k})\overset{k}{\to}0, \ \forall\ x\in B_\delta^c(y_0).
\end{equation}
\end{lem}

\noindent\textbf{Proof.} By
(\ref{eq2.18}) and (\ref{3:eq2.15}), we see that $w_k$ satisfies
\begin{equation}\label{3:eq2.19}
-\Delta w_k(x)+\varepsilon _k^2 (|\varepsilon _kx+\varepsilon _k
y_{\varepsilon _{k}} |-A)^2w_k(x)=\mu_k\varepsilon _k^2 w_k(x)+a_k
w_k^3(x)\quad \text { in } \ \mathbb{R}^2,
\end{equation}
where $\mu_k\in \R^2$ is a Lagrange multiplier. Similarly
to the proof of Lemma \ref{le2.3} (iii),  we can prove that there exists a
subsequence of $\{w_k\}$,  still denoted by  $\{w_k\}$,  such that (\ref{3.7}) holds  and
 $w_k\xrightarrow{k} w_0$ strongly in $H^1(\mathbb{R}^2)$ for
some positive function $w_0$ satisfying
\begin{equation}\label{3:eq2.20}
-\Delta w_0(x)=-\beta ^2 w_0(x)+a^* w_0^{3}(x)\quad \text { \ in } \
\mathbb{R}^2,
\end{equation}
where $\beta >0$ is a positive constant. Hence, for an arbitrary
large number $\alpha>2$,
\begin{equation}\label{unif}
\int_{|x|\geq R}|w_k|^\alpha dx\rightarrow0 \quad \text{as}\quad
R\rightarrow\infty \quad \text{uniformly for large }  k.
\end{equation}
Note from (\ref{3:eq2.19}) that $-\Delta w_k-c(x)w_k\leq 0$, where
$c(x)=a_k w_k^2(x)$. By applying De Giorgi-Nash-Moser theory, we
thus have $$\max_{B_1(\xi)} w_k\leq
C\Big(\int_{B_2(\xi)}|w_k|^\alpha dx\Big)^\frac{1}{\alpha},$$ where
$\xi$ is an arbitrary point in $\mathbb{R}^2$, and $C$ is a constant
depending only on the bound of $\|w_k\|_{L^\alpha(B_2(\xi))}$. We
hence deduce from (\ref{unif}) that
\begin{equation}
w_k(x)\rightarrow 0 \quad \text{as}\quad |x|\rightarrow\infty \quad
\text{uniformly in } k.
\end{equation}
Since $w_k$ satisfies (\ref{3:eq2.19}), one can use  the comparison
principle as in \cite{KW} to compare $w_k$ with
$Ce^{-\frac{\beta}{2}|x|}$, which then shows that there exists a
large constant $R>0$, independent of $k$, such that
\begin{equation}\label{exp}
w_k(x)\leq Ce^{-\frac{\beta}{2}|x|} \quad \text{for} \quad |x|>R
\quad \text{as}\quad k\rightarrow\infty.
\end{equation}

For any $x\in B_\delta^c(y_0)$, it then follows from (\ref{3.7}) that
$$\frac{|x-z_k|}{\eps_k}\geq\frac{1}{2}\frac{|x-y_0|}{\eps_k}\geq \frac{\delta}{2\eps_k}\overset{k}{\to}+\infty,$$
which, together with (\ref{exp}), yields that

$$u_k(x)=\frac{1}{\varepsilon_k}w_k(\frac{x-z_k}{\varepsilon_k})\leq \frac{1}{\varepsilon_k}e^{-\frac{\beta\delta}{2\varepsilon_k}}\overset{k}{\to}0,\ \forall \ x\in B_\delta^c(y_0),$$
i.e., (\ref{3.8}) holds.\qed\\

Motivated by \cite{GS,Wang},  we are now
ready to prove  Theorem \ref{thmA}.

\noindent\textbf{Proof of Theorem \ref{thmA}.} We still set $\varepsilon
_k=(a^*-a_k)^{\frac{1}{4 }}>0$, where $a_k \nearrow a^*$,  and  $u_k(x):=u_{a_k}(x)$ is a non-negative
minimizer   of (\ref{eq1.3}).  We start  the proof by establishing  first  the detailed concentration behavior of $u_k$.

 Let
$\bar{z}_k$ be any local maximum point of $u_k$. It then follows from
(\ref{3:eq2.19})   that \begin{equation}\label{3.14}u_k(\bar{z}_k)\geq
(\frac{-\mu_k}{a_k})^\frac{1}{2}\geq C\varepsilon_k^{-1}.\end{equation}
This
estimate and (\ref{3.8}) imply that, by passing to subsequence,
\begin{equation}\label{3.12}\bar z_k\xrightarrow{k}y_0\in\R^2\text{ with } |y_0|=A.\end{equation}
Set
\begin{equation}\label{3:W}\bar{w}_k=\varepsilon_ku_k(\varepsilon_k x+\bar{z}_k).\end{equation}
It then follows from (\ref{3:eq2.19}) that
\begin{equation}\label{3:222}
-\Delta \bar{w}_k(x)+\varepsilon_k^2 (|\varepsilon_k
x+\bar{z}_k|-A)^2\bar{w}_k(x)=\mu_k\varepsilon_k^2 \bar{w}_k(x)+a _k
\bar{w}_k^{3}(x) \ \text { in } \ \mathbb{R}^2.
\end{equation}
We claim that $\bar{w}_k$ satisfies (\ref{3:eq2.16}) for some
positive constants $R_0$ and $ \eta$. To prove this claim, we first
show that $\{\frac{\bar z_k-z_k}{\eps_k}\}$ is bounded  in
$\R^2$. Otherwise, if $\big|\frac{\bar z_k-z_k}{\eps_k}\big|\to
\infty$ as $k\to \infty$, it then follows  from the exponential
decay (\ref{exp}) that
\begin{equation*}
u_k(\bar z_k)=\frac{1}{\eps _k}w_k(\frac{\bar z_k-z_k}{\eps_k})\leq
\frac{C}{\eps _k}e^{-\frac{\beta}{2}\big|\frac{\bar
z_k-z_k}{\eps_k}\big|}=o( \eps_k^{-1}) \quad \text{as}\quad k\to
\infty\,,
\end{equation*}
which however contradicts (\ref{3.14}). Therefore, there exists a constant $R_1>0$,
independent of $k$, such that $\big|\frac{\bar
z_k-z_k}{\eps_k}\big|<\frac{R_1}{2}$. Note from (\ref{3:W}) and
(\ref{3:eq2.15}) that $$\bar w_k(x)=w_k(x+\frac{\bar
z_k-z_k}{\eps_k})\,.$$ Since $w_k$ satisfies (\ref{3:eq2.16}), we
then obtain from the above  that
\begin{equation}\label{3:est}
\lim_{k\to \infty}\int_{B_{R_0+R_1}(0)}|\bar w_k|^2dx=\lim_{k\to
\infty}\int_{B_{R_0+R_1}(\frac{\bar
z_k-z_k}{\eps_k})}|w_k|^2dx\geq\int_{B_{R_0}(0)}|w_k|^2dx\geq\eta>0\,,
\end{equation}
and the claim is therefore established.

As in the proof of Lemma \ref{le3.1} (see also Lemma \ref{le2.3} (iii)), one can further derive that there exists a
subsequence, still denoted by $\{\bar{w}_k\}$, of $\{\bar{w}_k\}$
such that
\begin{equation}\label{3.16}\bar{w}_k\xrightarrow{k} \bar{w}_0\text{ strongly in }H^1(\mathbb{R}^2),\end{equation}
for some nonnegative function $\bar{w}_0\geq 0$, where $\bar w_0$
satisfies (\ref{3:eq2.20}) for some  constant $\beta >0$. Note from
(\ref{3:est})  that $\bar w_0\not\equiv0$.  Thus, the strong maximum
principle yields that $\bar{w}_0(x)>0$ in $\R^2$. Since the origin
is a critical point of $\bar{w}_k$ for all $k>0$, it is also a
critical point of $\bar{w}_0$.  We therefore conclude from the
uniqueness (up to translations) of positive radial solutions for
(\ref{Kwong}) that   $\bar{w}_0$ is spherically symmetric about the
origin, and
\begin{equation}\label{eq2.350}
\bar w_0=\frac{\beta}{\|Q\|_2}Q(\beta|x|) \quad \text{for \, some }\
\beta >0.
\end{equation}

One can deduce from the above that $\bar{w}_k\geq
(\frac{\beta^2}{2a^*})^\frac{1}{2}$ at each local maximum point.
Since $\bar{w}_k$ decays to zero uniformly in $k$ as
$|x|\rightarrow\infty$, all   local maximum points of $\bar{w}_k$
stay in a finite ball in $\mathbb{R}^2$. Since $\bar{w}_k\rightarrow
\bar{w}_0$ in $C^2_{loc}(\mathbb{R}^2)$ and the origin is the only
critical point of $\bar{w}_0$, all local maximum points must
approach the origin and hence stay in a small ball $B_\epsilon(0)$
as $k\rightarrow\infty$. One can take $\epsilon$ small enough such
that $\bar{w}''_0(r)<0$ for $0\leq r\leq \epsilon$. It then follows
from Lemma 4.2 in \cite{NT} that for large $k$, $\bar{w}_k$ has no
critical points other than the origin. This  gives the uniqueness of
local maximum points
for $\bar{w}_k(x)$, which therefore implies that there exists a subsequence of $u_k$ concentrating at a {\em unique} global minimum point of potential $V(x)=(|x|-A)^2$.\\

To complete the proof of Theorem \ref{thmA}, we need to compute the exact  value of $\beta$ defined in (\ref{eq2.350}).  From (\ref{3:W}), we have
\begin{align}\label{3:eqa}
e(a_k) = E_{a_k}(u_k) & = \frac 1 {\eps_k^2} \Big[  \inte |\nabla
\bar{w}_k (x)|^2 dx - \frac {a^*} 2 \inte \bar{w}_k ^4(x) dx
\Big]\nonumber
\\ & \quad + \frac{\eps_k^2}2 \inte \bar{w}^4_k (x) dx + \inte
\big(|\eps_k x+\bar{z}_k |-|y_0|\big)^2 \bar{w}_k^2 (x) dx\,,
\end{align}
where $\bar{z}_k$ is the unique global maximum  point of $u_k$, and
$\bar{z}_k \to y_0\in \R^2$ as $k\to\infty$ for some    $|y_0|=A>0$.
The term in square brackets is non-negative and can be dropped for
the lower bound of $e(a_k)$. The $L^4(\R^2)$ norm of $\bar w_k $
converges to that of   $\bar w_0 $ as $k\to\infty$.

To estimate the last term of (\ref{3:eqa}), we   claim that
$\big\{\frac{|\bar{z}_k|-|y_0|}{\varepsilon_k}\big\}$ is bounded uniformly for
$k\to\infty$. Indeed, if there exists a subsequence of $\{a_k\}$, still denoted
by $\{a_k\}$, such that
$\big|\frac{|\bar{z}_k|-|y_0|}{\varepsilon_k}\big|\rightarrow\infty$
as $k\to\infty$, then for an arbitrary constant $C>0$,
 \begin{equation*}
\begin{split}
 &\lim_{k \to \infty} \eps_k^{-2} \inte \big(|\eps_k x+\bar z_k
|-|y_0|\big)^2 \bar w_k^2 (x) dx =\lim_{k \to \infty}\inte
\Big(|x+\frac{\bar z_k}{\eps_k}
|-\frac{|y_0|}{\varepsilon_k}\Big)^2\bar w_k^2 (x) dx \geq C\,.
\end{split}
\end{equation*}
This estimate and (\ref{3:eqa}) then imply that $$e(a_k)\geq C
\eps_k^2=C (a^*-a_k)^\frac{1}{2}   $$ holds  for arbitrary constant
$C>0$, which however contradicts Theorem \ref{thmB}, and the claim
is therefore true.

We now deduce from the above claim that there exists a subsequence
still denoted by $\{a_k\}$ such that
\begin{equation}\label{eq2.46a}
\frac{|\bar{z}_k|-|y_0|}{\varepsilon_k} \to C_0 \quad \text{as}\quad
k\rightarrow\infty \end{equation} for some constant $C_0$. Since $Q$
is a radial decreasing function and decays exponentially as
$|x|\to\infty$, we then deduce from (\ref{eq2.350}) that
\begin{align}\label{str}
 &\lim_{k \to \infty} \frac{1}{\eps_k^{2}} \inte \big(|\eps_k x+\bar z_k
|-|y_0|\big)^2 \bar w_k^2 (x) dx \nonumber\\
&= \lim_{k \to \infty} \inte \Big(\frac{|\eps_k x+\bar z_k
|}{\eps_k}-\frac{|y_0|}{\eps_k}\Big)^2 \bar w_k^2 (x)dx\nonumber\\
 &=\lim_{k \to \infty} \inte \Big(\frac{|\eps_k x+\bar
z_k |-|\bar z_k|}{\eps_k}+\frac{|\bar z_k|-|y_0|}{\eps_k}\Big)^2
\bar w_k^2 (x)dx \nonumber\\
 &=\inte \Big(\frac{y_0\cdot
x}{|y_0|}+C_0\Big)^2 \bar w_0^2(x)dx\geq \inte  \frac{|y_0\cdot
x|^2}{A^2} \bar w_0^2(x)dx\,,
\end{align}
where the equality holds if and only if $C_0=0$. We hence infer from
(\ref{3:eqa}) and (\ref{str}) that
\begin{equation}\label{mb}
\begin{split} \lim_{k\to \infty}\frac{e(a_k)}{(a^*-a_k)^{1/2}} &\geq
\frac 12 \|\bar w_0\|_4^4 + \frac{1}{A^2}\inte |y_0\cdot
x|^2  \bar w_0^2(x)dx \\
&= \frac{1}{a^*} \Big( \beta^2 + \frac{1}{A^2\beta^2} \inte
|y_0\cdot x |^2Q^2(x)dx \Big)\,,
\end{split}
\end{equation}
where  (\ref{1:id})  is used in the equality. Thus, taking the
infimum over $\beta>0$ yields that\begin{equation}\label{lim}
\lim_{k\to \infty}\frac{e(a_k)}{(a^*-a_k)^{1/2}} \geq\frac{2
}{a^*}\Big( \frac{\inte |y_0\cdot
x|^2Q^2(x)dx}{A^2}\Big)^\frac{1}{2}\,,
\end{equation} where the equality is achieved at  $$\beta= \lam _0:=\Big(
\frac{\inte |y_0\cdot x|^2Q^2(x)dx}{A^2}\Big)^\frac{1}{4}.$$

We finally note that the limit in (\ref{lim}) actually exists, and
it is equal to the right hand side of (\ref{lim}). To see this, one
simply takes
$$u (x)=\frac{\beta }{\eps \|Q\|_2}Q\Big(\frac{\beta |x-y_0|}{\eps}\Big) $$
as a trial function for $E_{a}(\cdot)$ and minimizes over $\beta >0$. The
result is that
\begin{equation}\label{lim2}
\lim_{a_k\nearrow a^*}\frac{e(a_k)}{(a^*-a_k)^{1/2}} =\frac{2}{a^*}\Big(
\frac{\inte |y_0\cdot x|^2Q^2(x)dx}{A^2}\Big)^\frac{1}{2}\,.
\end{equation}
The equality (\ref{lim2}) gives us two conclusions. Firstly, $\beta$
is unique, which is independent of the choice of the subsequence,
and equal to the expression minimizing (\ref{mb}), i.e., $\beta=\lam
_0$. Secondly, (\ref{str}) is indeed an equality, and thus $C_0=0$, i.e., (\ref{1.10}) holds. Moreover, from (\ref{3.12}), (\ref{3.16}) and (\ref{eq2.350}), we see that
\[
\bar{w}_k(x)=\varepsilon_k u (\varepsilon_k  x+\bar{z}_k)\xrightarrow{k} \frac{\lam
_0}{\|Q\|_2}Q(\lam _0|x|) \text{ strongly in }H^1(\R^2),
\]
where $\bar z_k$ is the unique maximum point of $u_k$ and  $\bar{z}_k\xrightarrow{k}y_0$  for some $y_0\in \R^2$
satisfying $|y_0|=A>0$. This completes the proof of Theorem
\ref{thmA}. \qed

Theorem \ref{thmA} gives a detailed description on the concentration behavior of  the minimizers of $e(a)$ when $a$ is close to $a^*$, upon which the phenomena of symmetry breaking of the the minimizers of $e(a)$ can be demonstrated in Corollary \ref{coro1} which we are going to prove now. With these conclusions we see that, when $a$ increases from $0$ to $a^*$, the minimizers of GP energy $e(a)$ have essentially different properties. The GP energy $e(a)$ has a unique non-negative minimizer which is radially symmetric if $a>0$ is small, but $e(a)$ has infinity many minimizers which are non-radially symmetric if $a$ approaches $a^*$.

{\bf Proof of Corollary \ref{coro1}.} {\bf (i):} By the proof of Theorem \ref{u:th}(ii) in the Appendix, we know that there exists a positive constant $a_{**}$ such that $e(a)$ has a unique non-negative minimizer for each $a\in [0,a_{**})$. Moreover, it is not difficult to see that this unique minimizer must be radially symmetric since $V(x)$ given by (\ref{1:V}) is radially symmetric and any rotation of the minimizer is still a minimizer of $e(a)$.

{\bf (ii):} It follows from Theorem \ref{thmA} that there exists $a_*>0$ such that, for $a\in (a_*,a^*)$, any non-negative minimizer of $e(a)$ has a unique maximum point $x_a$, and $\{x_a\}_a$ tends to a point on the  circular bottom of $V(x)$: $\{x\in\R^2: |x|=A\}$  as $a\nearrow a^*$. That is, when $a\in (a_*,a^*)$ each minimizer of $e(a)$ concentrates around a point on the circular bottom of $V(x)$, which then implies that all minimizers of $e(a)$ cannot be radially symmetric for $a\in (a_*,a^*)$. Moreover, since $V(x)$ is radially symmetric,  by rotation we can find that there are infinitely many minimizers of $e(a)$ for any $a\in (a_*,a^*)$. \qed


\small
\section{Appendix: Proof of Theorem \ref{u:th}(ii)}
In this appendix, we always assume that
\begin{equation}\label{4.1}0\leq V(x)\in L^\infty_{\rm
loc}(\R^2)\,,\quad\lim_{|x|\to\infty} V(x) = \infty\quad\text{and}\quad\inf_{x\in
\R^2} V(x) =0\,.\end{equation}
The properties of the Schr\"odinger operator $-\Delta+V(x)$ with $V(x)\in L^\infty(\mathbb{R}^2)$ are well known, see e.g. \cite{stuart}, but we could not find a reference for that of $V(x)$ satisfying (\ref{4.1}) although we guess  it  should exist somewhere. For the sake of completeness, we begin this appendix by  giving  some properties  of the Schr\"odinger operator $-\Delta+V(x)$ under conditions (\ref{4.1}),  which are required in proving the uniqueness of non-negative minimizers.  Before going to the properties of $-\Delta+V(x)$, we recall the following  embedding lemma, which can be found in \cite[Theorem  XIII.67]{RS}, \cite[Lemma 2.1]{DL} or \cite[Lemma 2.1]{Bao2}, etc.
\begin{lem}\label{1:lem1}
Suppose $V(x)$ satisfies   (\ref{4.1}).
Then, the embedding $\h\hookrightarrow
L^{q}(\R^2)$ is compact for all $q\in[2,\infty)$.
\end{lem}\qed

Define
\begin{equation}\label{u:mu1}\mu_1=\inf\Big\{\inte |\nabla u|^2+V(x)u^2dx:\  u\in\h \text{ and }\inte u^2dx=1\Big\}\,.\end{equation}
Usually, $\mu_1$ is called the first eigenvalue of $-\Delta+V(x)$ in $\h$. Moreover, by
Lemma \ref{1:lem1}, it is not difficult to know that $\mu_1$ is simple and can be attained by a positive function $\phi_1\in\h$. $\phi_1>0$ is then called the first eigenfunction of $-\Delta+V(x)$ in $\h$. We now define
\begin{equation}\label{u:mu2}\mu_2=\inf\Big\{\inte |\nabla u|^2+V(x)u^2dx: u\in Z\text{ and }\inte u^2dx=1\Big\}\,,\end{equation}
 where $$Z= \text{span}\{\phi_1\}^\perp=\Big\{u: u\in\h,\inte u \phi_1dx=0\Big\}\,.$$
 It is known that $\mu_2>\mu_1$ and
 \begin{equation}\label{u:sep}\h= \text{span}\{\phi_1\}\oplus Z\,.\end{equation}
Then, we have the following lemma.

\begin{lem}\label{u:lem2} Under the assumption of   (\ref{4.1}), we have
\begin{enumerate}
   \item [\rm(i)]  ker$\big(-\Delta+V(x)-\mu_1\big)=\text{span}\{\phi_1\}$;
   \item [\rm(ii)] $\phi_1\not\in \big(-\Delta+V(x)-\mu_1\big)Z$;
   \item [\rm(iii)] Im $\big(-\Delta+V(x)-\mu_1\big)$=$\big(-\Delta+V(x)-\mu_1\big)Z$ is closed in $\h^*$;
   \item[\rm(iv)] codim Im $\big(-\Delta+V(x)-\mu_1\big)$=1 .
 \end{enumerate}
 where $\h^*$ denotes the dual space of $\h$.
\end{lem}
\noindent\textbf{proof:}
  Since the proofs of (i) and (ii) are standard, we only prove (iii) and (iv).

{ \bf(iii):}  For any $z\in Z$, and $\psi \in \h$, we obtain from (\ref{u:mu1}) and (\ref{u:mu2}) that
  \begin{equation}\label{u:h1}
  \begin{split}
&\Big| \big\langle (-\Delta+V(x)-\mu_1)z, \psi\big\rangle_{\h^*,\h}\Big|=\Big|\inte \nabla z \nabla \psi+V(x)z \psi dx-\mu_1\inte z\psi dx\Big|\\
&\leq \|z\|_{\h}\|\psi\|_{\h}+\mu_1\|z\|_2\|\psi\|_2\leq(1+\sqrt\frac{\mu_1}{\mu_2})\|z\|_{\h}\|\psi\|_{\h}\,.
\end{split}
 \end{equation}
 On the other hand, it follows from  (\ref{u:mu1}) that
 \begin{equation}\label{u:h2}
  \begin{split}
 \big\langle (-\Delta+V(x)-\mu_1)z, z\big\rangle_{\h^*,\h}&=\inte |\nabla z|^2 +V(x)z^2 dx-\mu_1\inte z^2dx\\
 &\geq \big(1-\frac{\mu_1}{\mu_2}\big)\|z\|_{\h}^2\,.
\end{split}
 \end{equation}
It then follows  from  (\ref{u:h1}) and (\ref{u:h2})    that
 $$\big(1-\frac{\mu_1}{\mu_2}\big)\|z\|_{\h}\leq \|(-\Delta+V(x)-\mu_1)z\|_{\h^*}\leq (1+\sqrt\frac{\mu_1}{\mu_2})\|z\|_{\h}\quad \text{for all}\quad z\in Z,$$
which  proves (iii).

  {\bf(iv):} For any $f\in \h^*$, by Riesz representation theorem, there exists $u\in \h$ such that
 $$f(\psi)=\langle u,\psi\rangle_{\h,\h}\quad \text{for all}\quad \psi\in\h\,,$$
 i.e.,
 $$f=\big(-\Delta+V(x)\big)u\,.$$
 In view of  (\ref{u:sep}), we can write $u=c\phi_1+z_1$ for some  $z_1\in Z$. Thus,
 \begin{equation}\label{u:f}f=\big(-\Delta+V(x)\big)u=c\mu_1\phi_1+\big(-\Delta+V(x)\big)z_1 \quad \text{in}\quad \h^*\,.\end{equation}

Note from  (\ref{u:h1}) and (\ref{u:h2}) that   $\big\langle(-\Delta+V(x)-\mu_1)\cdot, \cdot \big\rangle_{\h^*,\h}$ is an inner product on $Z$, and $\big(-\Delta+V(x)\big)z_1\in Z^*$, where $Z^*$ denotes the dual space of $Z$. Hence, by applying Riesz representation theorem again, there exists $z_2\in Z$ such that
 \begin{equation}\label{u:allz}\big\langle(-\Delta+V(x))z_1, z\big\rangle_{\h^*,\h}=\big\langle(-\Delta+V(x)-\mu_1)z_2, z\big\rangle_{\h^*,\h}\quad \text{for all}\quad z\in Z\,.\end{equation}
 Moreover, since
 \begin{equation*}\begin{split}&\big\langle(-\Delta+V(x))z_1, \phi_1\big\rangle_{\h^*,\h}=\big\langle z_1, (-\Delta+V(x))\phi_1\big\rangle_{\h^*,\h}=\mu_1\inte z_1\phi_1=0\,,\\
& \big\langle(-\Delta+V(x)-\mu_1)z_2, \phi_1\big\rangle_{\h^*,\h}=\big\langle z_2, (-\Delta+V(x)-\mu_1)\phi_1\big\rangle_{\h^*,\h}=0,
 \end{split}\end{equation*}
We  obtain from (\ref{u:allz}) that
 \begin{equation*}\big\langle(-\Delta+V(x))z_1, \psi\big\rangle_{\h^*,\h}=\big\langle(-\Delta+V(x)-\mu_1)z_2, \psi\big\rangle_{\h^*,\h}\quad \text{for all}\quad \psi\in \h\,,\end{equation*}
 equivalently,
 $$\big(-\Delta+V(x)\big)z_1=\big(-\Delta+V(x)-\mu_1\big)z_2 \quad \text{in}\quad \h^*\quad \text{for some} \quad z_2\in Z\,.$$
  It then follows from (\ref{u:f}) that
  $$f=\big(-\Delta+V(x)\big)u=c\mu_1\phi_1+\big(-\Delta+V(x)-\mu_1\big)z_2 \quad \text{in}\quad \h^*\,.$$
 Then, by (ii), we   have
 $$\h^*=\text{span}\{\phi_1\}\oplus (-\Delta+V(x)-\mu_1)Z, $$
which implies  (iv).\qed\\

Stimulated  by Theorem 3.2 in \cite{CRa},  we have the following  lemma.
 \begin{lem}\label{u:le3}
 Define the following $C^1$ functional  $F: \h\times \R^2\mapsto \h^*$
 \begin{equation}\label{u:eq1}F(u,\mu, a)=\big(-\Delta +V(x)-\mu\big)u-a u^3\,.\end{equation}
Then, there exist $\delta>0$ and  a  unique function $\big(u(a),\mu(a)\big)\in C^1\big(B_\delta(0); B_\delta(\mu_1,\phi_1)\big)$ such that
\begin{equation} \label{u:F}
\begin{cases}
\mu(0)=\mu_1,\quad  u(0)=\phi_1;\\
F\big(u(a), \mu(a), a\big)=0;\\
\|u(a)\|_2^2=1\,.
\end{cases}
\end{equation}

 \end{lem}
 \noindent{\bf Proof.}  Let $g: Z\times\R^3\mapsto \h^*$ be defined  by
 $$g(z, \tau, s,a):=F\big((1+s)\phi_1+z,\mu_1+\tau, a\big)\,.$$
Then  $g\in C^1( Z\times\R^3, \h^*)$ and
\begin{equation}
 \begin{split}
 &g(0,0,0,0)=F(\phi_1, \mu_1,0)=0\,,\\
&g_{s}(0,0,0,0)=F_u(\phi_1, \mu_1, 0)\phi_1 =\big(-\Delta+V(x)-\mu_1\big)\phi_1=0\,.
 \end{split}
 \end{equation}
 Moreover, for any $(\hat z, \hat \tau)\in Z\times \R$ we have
 \begin{equation}
 \begin{split}
 g_{(z,\tau)}(0,0,0,0)(\hat z, \hat \tau)&=F_u(\phi_1, \mu_1, 0)\hat z +F_\mu(\phi_1, \mu_1, 0)\hat \tau\\&=\big(-\Delta+V(x)-\mu_1\big)\hat z-\hat\tau \phi_1\,.
 \end{split}
 \end{equation}
 Then, by Lemma \ref{u:lem2}, $g_{(z,\tau)}(0,0,0,0): Z\times \R\mapsto \h^*$ is an isomorphism.
Therefore, it follows from  the implicit function theorem that there exist $\delta_1>0$ and a unique function $(z(s,a),\tau(s,a))\in C^1(B_{\delta_1}(0,0);B_{\delta_1}(0,0))$ such that
  \begin{equation}\label{u:g}
 \begin{cases}
 g(z(s,a),\tau(s,a), s,a)=F((1+s)\phi_1+z(s,a),\quad \mu_1+\tau(s,a),a)=0\,,\\
 z(0,0)=0,\tau(0,0)=0\,,\\
 z_s(0,0)=-g_{(z,\tau)}^{-1}(0,0,0,0)\cdot g_{s}(0,0,0,0)=0\,.
 \end{cases}
 \end{equation}

Now, let $$u(s,a)=(1+s)\phi_1+z(s,a)\,,\quad (s,a)\in B_{\delta_1}(0,0)\,,$$
 and define
 $$f(s,a)=\|u(s,a)\|_2^2=(1+s)^2+\inte z(s,a)^2dx\,,\quad (s,a)\in B_{\delta_1}(0,0)\,.$$
 By (\ref{u:g}), we have
 $$f(0,0)=1, \quad f_s(0,0)=2+2\inte z_s(0,0)z(0,0)dx=2\,.$$
 Then, by applying  implicit function theorem again, there exist $0<\delta<\delta_1$ and a unique function  $s=s(a)\in C^1\big(B_{\delta}(0);B_{\delta}(0)\big)$ such that
 $$f(s(a),a)=\|u(s(a),a)\|_2^2=f(0,0)=1\,,\quad a\in B_{\delta}(0) \,.$$
This and (\ref{u:g}) imply that for $a\in B_{\delta}(0)$, there exists a  unique  function:$$\Big(u(a):=u\big(s(a),a\big), \mu(a):=\mu_1+\tau \big(s(a),a\big)\Big)\in C^1\big(B_{\delta}(0);B_{\delta}(\phi_1,\mu_1)\big)$$
such that (\ref{u:F}) holds, and the proof is therefore complete.\qed\\


We are now ready  to prove the uniqueness of non-negative minimizers of (\ref{eq1.3}), that is, Theorem \ref{u:th} (ii).

\noindent\textbf{Proof of Theorem \ref{u:th} (ii):}
Let $u_a(x)>0$ be a non-negative minimizer of $e(a)$ with $a\in[0,a^*)$.
It is easy to see that
\begin{equation}\label{u:ea}e(0)=\mu_1 \quad \text{and}\quad e(a)\leq e(0)=\mu_1\,,\end{equation}
where $\mu_1$ is the first eigenvalue of $-\Delta+V(x)$.
We  first prove that
\begin{equation}\label{u:co}e(a)\in C([0,a^*),\R^+)\,.\end{equation}
Indeed, for any $a_0\in [0,a^*)$, it follows from Gagliardo-Nirenberg inequality
(\ref{GNineq}) that
\begin{equation}\label{u:bound}\inte u_a^4dx\leq \frac{2e(a)}{a^*-a}\leq\frac{2\mu_1}{a^*-a}\leq \frac{4\mu_1}{a^*-a_0} \quad \text{for}\quad 0\leq a\leq \frac{a^*+a_0}{2}\,.\end{equation}
Thus,
$$e(a_0)\leq E_{a_0}(u_a)=E_{a}(u_a)+\frac{a-a_0}{2}\inte u_a^4dx=e(a)+O(|a-a_0|)\quad \text{as}\quad a\to a_0\,,$$
$$e(a_0)= E_{a_0}(u_{a_0})=E_{a}(u_{a_0})+\frac{a-a_0}{2}\inte u_{a_0}^4dx\geq e(a)+O(|a-a_0|)\quad \text{as}\quad a\to a_0\,.$$
Hence,  $$\lim_{a\to a_0}e(a)=e(a_0)\,,$$
which  implies (\ref{u:co}).

Since $u_a$ is a  minimizer of (\ref{eq1.3}), it satisfies
the Euler-Lagrange equation
\begin{equation*}
-\Delta u_a(x)+V(x)u_a(x)-\mu_a u_a(x)-a u_a^3(x)=0\quad  \text{ in }\
\mathbb{R}^2,
\end{equation*}
where $\mu_a\in \mathbb{R}$ is a suitable Lagrange multiplier, i.e.,
\begin{equation}\label{4.17}
F(u_a,\mu_a,a)=0, \quad \text{where $F(\cdot)$ is defined  by (\ref{u:eq1})}.
\end{equation}
Since
\begin{equation*}
\mu_a=e(a)-\frac{a}{2}\int_{\mathbb{R}^2}|u_a|^4dx,
\end{equation*}
it then follows from (\ref{u:ea})-(\ref{u:bound}) that there exists $a_1>0$ small such that
\begin{equation}\label{u:mua}|\mu_a-\mu_1|\leq |e(a)-\mu_1|+\frac{a}{2}\int_{\mathbb{R}^2}|u_a|^4dx\leq \delta \quad \text{for}\quad 0\leq a<a_1\,,\end{equation}
where $\delta>0$ is as in Lemma \ref{u:le3}.
On the other hand, since
$$E_0(u_a)=e(a)+\frac{a}{2}\int_{\mathbb{R}^2}|u_a|^4dx\to e(0)=\mu_1\quad \text{as}\quad a\searrow 0\,,$$
i.e., $\{u_a\geq0\}$  is a minimizing sequence of $e(0)=\mu_1$ as $a\searrow 0$. Noting that   $\mu_1$ is a simple eigenvalue, we can easily deduce from  Lemma \ref{1:lem1} that $$u_a\to \phi_1\quad \text{in}\quad \h\quad \text{for all}\quad a\searrow 0\,.$$
This implies that there exists $a_2>0$ such that \begin{equation}\label{u:clo}\|u_a-\phi_1\|_{\h}<\delta\quad \text{for}\quad 0\leq a<a_2\,.\end{equation}
Then using  (\ref{4.17})-(\ref{u:clo}) and   Lemma \ref{u:le3}, we obtain  that
$$\mu_a= \mu(a);\quad u_a=u(a)\quad \text{for}\quad 0\leq a<\min\{a_1,a_2\}\,,$$
i.e., $e(a)$ has a  unique non-negative minimizer $u(a)$  if $a>0$ is small.
\qed\\

\vspace {.95cm}

\normalsize

\noindent {\bf Acknowledgements:} The authors would like  to thank Professor  Robert
Seiringer for his fruitful discussions. Thanks also go to Professor Charles Stuart for helping us to prepare the final version of the paper.  This research was  supported
    by National Natural Science Foundation of China (11322104, 11171339, 11271360) and National Center for Mathematics and Interdisciplinary Sciences.

\ \\
\ \\
\ \\
{\it
\noindent Wuhan Institute of Physics and Mathematics,\\
Chinese Academy of Sciences,\\
P.O. Box 71010, Wuhan 430071,\\
 P. R. China\\}
\noindent email: yjguo@wipm.ac.cn\\
email: zeng6656078@126.com\\
email: hszhou@wipm.ac.cn
\end{document}